# TESTING FOR CHANGE POINTS IN TIME SERIES MODELS AND LIMITING THEOREMS FOR NED SEQUENCES[1]

BY SHIQING LING

*Hong Kong University of Science and Technology*

This paper first establishes a strong law of large numbers and a strong invariance principle for forward and backward sums of near-epoch dependent sequences. Using these limiting theorems, we develop a general asymptotic theory on the Wald test for change points in a general class of time series models under the no change-point hypothesis. As an application, we verify our assumptions for the long-memory fractional ARIMA model.

**1. Introduction.** Testing on structural change problems has been an important issue in statistics. The earliest references go back to Chow [9] and Quandt [33]. Chow's test is to assume that the time of structural change is known a priori, and the critical values for the $\chi^2$ distribution can be simply used. Quandt's test is to take the largest Chow test statistic over all possible times of the structural change. Quandt's test appears to be more reasonable in practice because it does not need to assume the time of structural change a priori. However, its critical values are hard to obtain even approximately due to singular behavior near the end points. One method is to restrict the change-point interval $(0,1)$ to $[\tau_1, \tau_2]$ with $0 < \tau_1 < \tau_2 < 1$; see [[2], [4], [15], [16], [22]]. Another important method is to normalize the Quandt-type test. This type of test statistic has a Darling–Erdös-type limit and its critical values are easily obtained. This method was developed by Yao and Davis [38] for i.i.d. normal data, and was extended by Horváth [17] for general i.i.d. data and Horváth [18] for linear regression models. However, when using this method for time series models, we encounter some great challenges.

To understand these, we look at the AR(1) model, $y_t = \phi y_{t-1} + \varepsilon_t$, where $|\phi| < 1$ and $\{\varepsilon_t\}$ are independent and identically distributed (i.i.d.) errors.

Received December 2004; revised September 2006.
[1]Supported by Hong Kong RGC (CERG No. HKUST4765/03H, No. HKUST602205/05P and No. HKUST6428/06H).

First, we need to obtain the rate of uniform convergence of the partial sample information matrices based on $\{y_1, \ldots, y_k\}$ and $\{y_{k+1}, \ldots, y_n\}$, respectively; that is, for some $\delta > 0$, we need to establish

(a) $\displaystyle\max_{g_n \leq k \leq n} \left| \frac{1}{k^{1-\delta}} \sum_{t=1}^{k} X_t \right| = o_p(1)$ and (b) $\displaystyle\max_{g_n \leq \tilde{k} < n} \left| \frac{1}{\tilde{k}^{1-\delta}} \sum_{t=k+1}^{n} X_t \right| = o_p(1)$,

as $n \to \infty$, where $\tilde{k} = n - k$, $g_n = \log\log\log(\max\{e^e, n\})$ and $X_t = y_{t-1}^2 - Ey_{t-1}^2$ (see Lemma 6.2). Under the strong mixing condition with $Ey_t^4 < \infty$, Davis, Huang and Yao [12] first established that

(1.1) $$\frac{1}{k^{1-\delta}} \sum_{t=1}^{k} X_t = o(1) \quad \text{a.s.,}$$

using the strong invariance principle in Kuelbs and Philipp [23], and then used (1.1) to obtain (a). We note that the ergodic theorem only ensures that $\sum_{t=1}^{k} X_t / k = o(1)$ a.s., which cannot be used for (a), and hence, (1.1) in [12] is novel. Since $\{y_t\}$ is strictly stationary, (b) is equivalent, for any $\epsilon > 0$, to

$$P\left( \max_{g_n \leq \tilde{k} < n} \frac{1}{\tilde{k}^{1-\delta}} \left| \sum_{t=-\tilde{k}+1+n}^{n} X_t \right| > \epsilon \right) = P\left( \max_{g_n \leq k < n} \frac{1}{k^{1-\delta}} \left| \sum_{t=-k}^{-1} X_t \right| > \epsilon \right) = o(1).$$

This is not equivalent to (a) if $\{y_t\}$ is not time-reversible. Except for Gaussian linear processes, very few time series have been shown to be time-reversible; see [8]. Thus, (1.1) cannot be used for (b), generally. To solve this problem, we need the following strong law of large numbers (SLLN):

(1.2) $$\frac{1}{k^{1-\delta}} \sum_{t=-k}^{-1} X_t = o(1) \quad \text{a.s.}$$

However, this has not been established in the literature.

Second, we need to approximate the score functions based on the subsamples $\{y_1, \ldots, y_k\}$ and $\{y_{k+1}, \ldots, y_n\}$ by i.i.d. normal random sequences $\{G_{1t} : t = 1, 2, \ldots\}$ and $\{G_{2t} : t = 1, 2, \ldots\}$, respectively, such that

(c) $\displaystyle\max_{g_n \leq k \leq n} k^\delta \left| \frac{1}{\sqrt{k}} \sum_{t=1}^{k} y_{t-1}\varepsilon_t - \frac{1}{\sqrt{k}} \sum_{t=1}^{k} G_{1t} \right| = o_p(1)$,

(d) $\displaystyle\max_{g_n \leq k < n} k^\delta \left| \frac{1}{\sqrt{k}} \sum_{t=-k}^{-1} y_{t-1}\varepsilon_t - \frac{1}{\sqrt{k}} \sum_{t=-k}^{-1} G_{2t} \right| = o_p(1)$,

for some $\delta > 0$. Davis, Huang and Yao [12] first used the result in Kuelbs and Philipp [23] to establish the strong invariance principle (SIP),

(1.3) $$\frac{1}{\sqrt{k}} \sum_{t=1}^{k} y_{t-1}\varepsilon_t = \frac{1}{\sqrt{k}} \sum_{t=1}^{k} G_{1t} + o(k^{-\delta}) \quad \text{a.s,}$$

with strong mixing $\{y_t\}$, for some $\delta > 0$, and then used (1.3) to prove (c). Similarly to (b), to prove (d), we need the backward SIP, that is, there is an i.i.d. normal random sequence $\{G_{2t}: t = 1, 2, \ldots\}$ such that

$$(1.4) \qquad \frac{1}{\sqrt{k}} \sum_{t=-k}^{-1} y_{t-1} \varepsilon_t = \frac{1}{\sqrt{k}} \sum_{t=1}^{k} G_{2t} + o(k^{-\delta}) \qquad \text{a.s.}$$

Again, there is not any result for (1.4) in the literature. The preceding difficulties are not only in Quandt-type tests but also in the estimating change-points as in [[3], [26], [31]]. This issue seems to be not well discussed in the literature.

This paper first establishes a new SLLN and a new SIP for the backward sums of near-epoch dependent (NED) sequences. The existing SLLN and SIP for the forward sums of random sequences related to (a) and (c), such as those in [34] and [14], require some mixing and high-order moment conditions, or do not have a rate of convergence (see also [25]). The mixing conditions are not always easy to verify. The high-order moment condition directly links to the restriction on the parameter space in some nonlinear time series models such as ARCH-type models. The weakest moment condition is in the ergodic theorem, but it does not have a rate of convergence. This paper next establishes a SLLN and a SIP with a rate of convergence for the forward sums of NED sequences under a weak moment condition and without a strong mixing assumption.

Our SLLNs and SIPs are given in Section 2. Using them, we study the Wald test for change-points in a class of time series models in Section 3. This is a general theory and can be used for many time series models. As an application, we verify our assumptions for long-memory FARIMA models in Section 4. The proofs are given in Sections 5–7. Throughout this paper, we use the following notation: $|A| = [\text{tr}(AA')]^{1/2}$ for a vector or matrix $A$ and $\|Z\|_p = (E|Z|^p)^{1/p}$ for a random vector or matrix $Z$ with its elements in $L^p$ space ($p \geq 1$). Finally, we refer to the related references [20] and [19] for Quandt-type tests with the long-memory time series, and to [24] for the sequential approach.

**2. Limiting theorems for NED sequences.** Let $\{\varepsilon_t\}$ be a series of independent random variables (or vectors) on the probability space $(\Omega, \mathcal{B}, P)$, $\mathcal{F}_t = \sigma\{\varepsilon_t, \varepsilon_{t-1}, \ldots\}$ and $X_t$ be a $\mathcal{F}_t$-measurable $m \times 1$ random vector for $t = 0, \pm 1, \ldots$. We first introduce the following definition.

DEFINITION 2.1. Let $\mathcal{F}_i(j)$ be the $\sigma$-field generated by $\{\varepsilon_j, \varepsilon_{j-1}, \ldots, \varepsilon_{j-i+1}\}$ with $i \geq 1$, and $\mathcal{F}_0(j) = \{\varnothing, \Omega\}$. $\{X_t\}$ is said to be $L^p(\nu)$ NED in terms of $\{\varepsilon_t\}$ if $\sup_{-\infty < t < \infty} \|X_t\|_p < \infty$ and $\sup_{-\infty < t < \infty} \|X_t - E[X_t|\mathcal{F}_k(t)]\|_p = O(k^{-\nu})$, where $p \geq 1$ and $\nu > 0$.

This notion of NED sequence extends a concept introduced in Billingsley [7]. Some different versions appear in [[30], [32], [36]]. This NED $\{X_t\}$ implies that it is mixingale, that is, $\sup_{-\infty<t<\infty}\|EX_t - E(X_t|\mathcal{F}_{t-k})\|_p = O(k^{-\nu})$. Our SLLN and SIP are as follows.

THEOREM 2.1. *Let $\{X_t : t = 0, \pm 1, \ldots\}$ be an $L^{1+\iota}(\nu)$ NED and mean zero sequence in terms of $\{\varepsilon_t\}$ with $\iota > 0$ and $\nu > 0$. Then there exists a constant $\delta > 0$ such that*

(a) $\dfrac{1}{k}\sum_{t=1}^{k} X_t = o\left(\dfrac{1}{k^\delta}\right)$ *a.s. and* (b) $\dfrac{1}{k}\sum_{t=-k}^{-1} X_t = o\left(\dfrac{1}{k^\delta}\right)$ *a.s.*

REMARK 2.1. The moment condition in Theorem 2.1 is only slightly stronger than that in the ergodic theorem for the forward sums. But our SLLN includes a rate of convergence, while the ergodic theorem does not. We guess that this is the weakest moment condition for the NED sequence if a rate of convergence is wanted. This rate is indispensable when we prove Lemmas 6.2–6.4. The independence of $\{\varepsilon_t\}$ can be replaced by some mixing conditions. If we allow $\iota \geq 1$ and $\nu \geq 0.5$, then a sharper rate of convergence may be obtained; see, for example, [14], page 41. If we assume $\iota = 1$ and use the moment bound of Ing and Wei [21], then a relationship between the rate of convergence and the series dependence can be given.

THEOREM 2.2. *Let $X_t$ be a martingale difference in terms of $\mathcal{F}_t$ with covariance matrix $\Omega$ and be $L^{2+\iota}(\nu)$ NED in terms of $\{\varepsilon_t\}$ with $\iota > 0$, where either $2\nu > 1$ or $2\nu = 1$ and there exist constants $\nu_1 > 0$ and $\iota_1 > 0$ with $2\nu_1 > 1$ such that*

(2.1) $$\sup_{-\infty<t<\infty} \|E[X_t|\mathcal{F}_{k+1}(t)] - E[X_t|\mathcal{F}_k(t)]\|_{2+\iota_1} = O(k^{-\nu_1}).$$

*Then, without changing its distribution, we can redefine $\{X_t\}$ on two richer probability spaces together with two sequences of i.i.d. $m \times 1$ normal vectors with mean zero and covariance matrix $\Omega$, $\{G_{1t} : t = 1, 2, \ldots\}$ and $\{G_{2t} : t = 1, 2, \ldots\}$, such that, for some constant $\delta > 0$, we have, respectively,*

(a) $\sum_{t=1}^{k} X_t = \sum_{t=1}^{k} G_{1t} + O(k^{1/2-\delta})$ *a.s. and*

(b) $\sum_{t=-k}^{-1} X_t = \sum_{t=1}^{k} G_{2t} + O(k^{1/2-\delta})$ *a.s.*

REMARK 2.2. The two richer probability spaces may be different, for which we refer to [6] and [13]. Theorems 2.1–2.2 do not require $\{X_t\}$ to be stationary and can be extended for triangular arrays as in [1].

**3. Testing change-points in time series models.** Assume that the time series $\{y_t : t = 0, \pm 1, \pm 2, \ldots\}$ is $\mathcal{F}_t$-measurable, strictly stationary and ergodic, and is generated by the model

$$y_t = f(\lambda, Y_{t-1}, \varepsilon_t), \tag{3.1}$$

where $f$ is a known function, $\lambda$ is an $m \times 1$ unknown parameter vector, $\{\varepsilon_t\}$ is i.i.d. and $Y_t = (y_t, y_{t-1}, \ldots)$. The structure of $\{y_t\}$ is characterized by $f$ and $\lambda$. This class of models (3.1) includes many time series models in the literature, such as ARMA, GARCH and random coefficient AR models. We assume that the parameter space $\Theta$ is a compact subset of $R^m$, and the true value of $\lambda$, denoted by $\lambda_0$, is an interior point in $\Theta$, where $R = (-\infty, \infty)$.

We denote the model (3.1) with the true parameter $\lambda_0$ by $M(\lambda_0)$. Let $y_1, \ldots, y_n$ be the observations. We consider the null and alternative hypotheses,

$$H_0 : \{y_1, \ldots, y_n\} \in M(\lambda_0) \quad \text{versus}$$

$$H_{1n}(k) : \{y_1, \ldots, y_k\} \in M(\lambda_0) \quad \text{and} \quad \{y_{k+1}, \ldots, y_n\} \in M(\lambda_{10})$$

$$\text{with } \lambda_0 \neq \lambda_{10} \text{ for some } k \in [1, n).$$

Here, $k = [n\tau]$ is called the change-point with $\tau \in (0,1)$, where $[x]$ is the integer part of $x$. Under $H_{1n}(k)$, we use the following objective functions (OF) to estimate $\lambda_0$ and $\lambda_{10}$, based on the sample $\{y_1, \ldots, y_n\}$ with initial value $Y_0$, respectively,

$$L_n(k, \lambda) = \sum_{t=1}^{k} l(\lambda, Y_t) \quad \text{and} \quad L_{1n}(k, \lambda_1) = \sum_{t=k+1}^{n} l(\lambda_1, Y_t), \tag{3.2}$$

where $l(\lambda, Y_t)$ is a measurable function in terms of $Y_t$ and is almost surely (a.s.) three times differentiable with respect to $\lambda$. The function $l(\lambda, Y_t)$ can be taken as those in LSE, MLE, quasi-MLE and $M$-estimators, among others. Let $l_t(\lambda) = l(\lambda, Y_t)$, $D_t(\lambda) = \partial l_t(\lambda)/\partial \lambda$ and $P_t(\lambda) = -\partial^2 l_t(\lambda)/\partial\lambda\,\partial\lambda'$. Denote $\Sigma = E[P_t(\lambda_0)]$ and $\Omega = E[D_t(\lambda_0)D_t'(\lambda_0)]$. Here and below, the expectation is with respect to the probability measure under the null hypothesis. We first give two sets of assumptions as follows.

ASSUMPTION 3.1. For some constant $\iota > 0$ and an open neighborhood $\Theta_0$ of $\lambda_0$:

(i) $E\sup_{\lambda \in \Theta} |l_t(\lambda)|^{1+\iota} < \infty$ and $E[l_t(\lambda)]$ has a unique maximum at $\lambda = \lambda_0$;

(ii) $D_t(\lambda_0)$ is an $\mathcal{F}_t$-measurable martingale difference with $\Omega > 0$;

(iii) $\Sigma > 0$ and $E\sup_{\lambda \in \Theta_0} |P_t(\lambda)|^{1+\iota} < \infty$;

(iv) $E\sup_{\lambda \in \Theta_0} |\partial p_{ijt}(\lambda)/\partial \lambda|^{0.5+\iota} < \infty$, where $p_{ijt}(\lambda)$ is the $(i,j)$th element of $P_t(\lambda)$.

ASSUMPTION 3.2. For some $\iota > 0$, $\nu_0 > 0$ and $\nu \geq 1/2$ and an open neighborhood $\Theta_0$ of $\lambda_0$:

(i) $\|\sup_{\lambda \in \Theta} |l_t(\lambda) - E[l_t(\lambda)|\mathcal{F}_k(t)]|\|_{1+\iota} = O(k^{-\nu_0})$;

(ii) $D_t(\lambda_0)$ is $L^{2+\iota}(\nu)$ NED in terms of $\{\varepsilon_t\}$ with $2\nu > 1$ or with $2\nu = 1$ and (2.1) being satisfied as $X_t = D_t(\lambda_0)$;

(iii) $\|\sup_{\lambda \in \Theta_0} |P_t(\lambda) - E[P_t(\lambda)|\mathcal{F}_k(t)]|\|_{1+\iota} = O(k^{-\nu_0})$.

When $\iota = 0$, Assumption 3.1(i)–(iii) is typical for estimating $\lambda_0$ in model (3.1). We need the $(1+\iota)$th finite moment here because the ergodic theorem cannot be used for $\hat{\lambda}_{1n}(k)$. Assumption 3.1(iv) is for the rate of uniform convergence in (6.2). Assumption 3.2 is a key to using Theorems 2.1–2.2 for Lemmas 6.1–6.2. In practice, $Y_0$ is usually replaced by some constants. Let $\widetilde{l}_t(\lambda)$, $\widetilde{D}_t(\lambda)$ and $\widetilde{P}_t(\lambda)$ be defined as $l_t(\lambda)$, $D_t(\lambda)$ and $P_t(\lambda)$, respectively, with initial values $y_t$ being zero or a constant for $t \leq 0$. Our initial condition is as follows.

ASSUMPTION 3.3. For some constants $\delta > 0$, $\nu_0 > 0$ and $\nu \geq 1/2$ and an open neighborhood $\Theta_0$ of $\lambda_0$:

(i) $E \sup_{\lambda \in \Theta} |l_t(\lambda) - \widetilde{l}_t(\lambda)| = O(t^{-\nu_0})$;

(ii) $\sup_{g_n \leq k \leq n} \{k^{-1/2+\delta} |\sum_{t=1}^{k} [D_t(\lambda_0) - \widetilde{D}_t(\lambda_0)]|\} = o_p(1)$ and $\sup_{g_n \leq n-k < n} \{(n-k)^{-1/2+\delta} |\sum_{t=k+1}^{n} [D_t(\lambda_0) - \widetilde{D}_t(\lambda_0)]|\} = o_p(1)$;

(iii) $E \sup_{\lambda \in \Theta_0} |P_t(\lambda) - \widetilde{P}_t(\lambda)| = O(t^{-\nu_0})$ and $\|D_t(\lambda_0) - \widetilde{D}_t(\lambda_0)\|_{1+\iota} = O(u_t)$, where $g_n = \log \log \log(\max\{e^e, n\})$ and $u_t = t^{-\nu} \log^q t$ for some $q > 0$.

It can be shown that Assumption 3.3(ii) holds if $2\nu > 1$ in the second part of Assumption 3.3(iii). The OFs in (3.2) are modified as

$$(3.3) \quad L_n(k, \lambda) = \sum_{t=1}^{k} \widetilde{l}_t(\lambda) \quad \text{and} \quad L_{1n}(k, \lambda_1) = \sum_{t=k+1}^{n} \widetilde{l}_t(\lambda_1).$$

Let $\hat{\lambda}_n(k)$ and $\hat{\lambda}_{1n}(k)$, respectively, be the maximizers of $L_n(k, \lambda)$ and $L_{1n}(k, \lambda_1)$ on $\Theta$ for each known $k$. The Wald test statistic evaluated at $[\hat{\lambda}_n(k), \hat{\lambda}_{1n}(k)]$ for testing $H_0$ against $H_{1n}(k)$ is defined as

$$W_n(k) = \frac{k(n-k)}{n^2} [\hat{\lambda}_n(k) - \hat{\lambda}_{1n}(k)]' [\hat{\Sigma}_n(k) \hat{\Omega}_n^{-1}(k) \hat{\Sigma}_n(k)] [\hat{\lambda}_n(k) - \hat{\lambda}_{1n}(k)],$$

where $\hat{\Sigma}_n(k) = \sum_{t=1}^{k} \widetilde{P}_t(\hat{\lambda}_n(k)) + \sum_{t=k+1}^{n} \widetilde{P}_t(\hat{\lambda}_{1n}(k))$ and

$$\hat{\Omega}_n(k) = \sum_{t=1}^{k} \widetilde{D}_t(\hat{\lambda}_n(k)) \widetilde{D}_t'(\hat{\lambda}_n(k)) + \sum_{t=k+1}^{n} \widetilde{D}_t(\hat{\lambda}_{1n}(k)) \widetilde{D}_t'(\hat{\lambda}_{1n}(k)).$$

When we test the null $H_0$ against $\bigcup_{k\in[1,n)} H_{1n}(k)$, a natural test statistic is $\max_{k\in[1,n)} W_n(k)$. However, this test statistic diverges to infinity; see [2]. We define the normalized Quandt-type Wald test statistic as

$$\widehat{W}_n(m) = \max_{m<k<n-m} \frac{W_n(k) - b_n(m)}{a_n(m)}, \tag{3.4}$$

where $a_n(m) = \sqrt{b_n(m)/(2\log\log n)}$, $b_n(m) = [2\log\log n + (m\log\log\log n)/2 - \log\Gamma(m/2)]^2/(2\log\log n)$ and $\Gamma(\cdot)$ is the gamma function. Our result for testing for a change-point in model (3.1) is as follows.

THEOREM 3.1. *If Assumptions* 3.1–3.3 *are satisfied, then under the null $H_0$, for any $x \in R$, $P[\widehat{W}_n(m) \le x] \to \exp(-2e^{-x/2})$ as $n \to \infty$.*

REMARK 3.1. Some weighted test statistics can be constructed along the lines of [11] where optimality of related tests is also discussed. Assumptions 3.1(i)–(iii) and 3.2–3.3 were verified by Ling [27] for the AR-GARCH model.

**4. Application to LM-FARIMA models.** The time series $\{y_t\}$ is said to follow a long-memory FARIMA$(p,d,q)$ model if

$$\phi(B)(1-B)^d y_t = \psi(B)\varepsilon_t, \tag{4.1}$$

where $\phi(B) = 1 - \sum_{i=1}^{p} \phi_i B^i$, $\psi(B) = 1 + \sum_{i=1}^{q} \psi_i B^i$, $B$ is the backward-shift operator, $d \in (0, 0.5)$ and $(1-B)^d = \sum_{k=0}^{\infty} c_k B^k$ with $c_k = (-d)(-d+1)\cdots(-d+k-1)/k!$, and $\{\varepsilon_t\}$ a sequence of i.i.d. white noise variables. Denote $\lambda = (d, \phi_1, \ldots, \phi_p, \psi_1, \ldots, \psi_q)'$. The parameter space $\Theta$ is a compact subset of $R^{p+q+1}$. Assume that the true parameter $\lambda_0$ of $\lambda$ is an interior point in $\Theta$ and, for each $\lambda \in \Theta$, it satisfies:

ASSUMPTION 4.1. $d \in (0, 0.5)$, $\phi(z) \ne 0$ and $\psi(z) \ne 0$ for all $z$ such that $|z| \le 1$, $\phi_p \ne 0$, $\psi_q \ne 0$, and $\phi(z)$ and $\psi(z)$ have no common root.

It is not hard to see that (4.1) is a special form of model (3.1). Following common practice, we use quasi-log-likelihood estimation for $\lambda_0$ and the OFs are

$$L_n(k,\lambda) = -\tfrac{1}{2}\sum_{t=1}^{k} \varepsilon_t^2(\lambda) \quad \text{and} \quad L_{1n}(k,\lambda_1) = -\tfrac{1}{2}\sum_{t=k+1}^{n} \varepsilon_t^2(\lambda_1), \tag{4.2}$$

where $\varepsilon_t(\lambda) = \psi^{-1}(B)\phi(B)(1-B)^d y_t$. In this case, we have

$$D_t(\lambda) = -\frac{\partial \varepsilon_t(\lambda)}{\partial \lambda}\varepsilon_t(\lambda) \quad \text{and} \quad P_t(\lambda) = \frac{\partial \varepsilon_t(\lambda)}{\partial \lambda}\frac{\partial \varepsilon_t(\lambda)}{\partial \lambda'} + \frac{\partial^2 \varepsilon_t(\lambda)}{\partial \lambda \partial \lambda'}\varepsilon_t(\lambda).$$

Let $\hat{\lambda}_n(k)$ and $\hat{\lambda}_{1n}(k)$ be the maximizers of $L_n(k,\lambda)$ and $L_{1n}(k,\lambda_1)$ on $\Theta$ for each $k$ with the initial values $y_t = 0$ for $t \le 0$. The result for model (4.1) is as follows.

THEOREM 4.1. *If Assumption 4.1 holds and $E|\varepsilon_t|^{2+\iota} < \infty$ for some $\iota > 0$, then under the null $H_0$, for any $x$, $P[\widehat{W}_n(p+q+1) \leq x] \to \exp(-2e^{-x/2})$ as $n \to \infty$.*

REMARK 4.1. For the linear processes with the long-memory parameter $H = 1/2 + d_0$, Beran and Terrin [5] and Horváth and Shao [20] proposed some tests for the change of $H$ in the frequency domain, but they did not verify the conditions for model (4.1) and assumed that $E|\varepsilon_t|^{4+\iota} < \infty$. See also [19]. As far as we know, our test statistic is new in the time domain and is also different from the tests in [5] and [20].

REMARK 4.2. To see the performance of the Wald test in finite samples, we examine a small simulation for the FARIMA$(0,d,0)$ model with $\varepsilon_t \sim N(0,1)$, using Fortran 77. Sample sizes $n = 250$ and $400$ are used. We first study the size, for which we take $d_0 = 0.1, 0.2, 0.3$ and $0.4$, and then the power, for which we take $d_0 = 0.1$ and $d_{10} = 0.2, 0.3, 0.4$ with the change-point $k = [0.5n]$ and $[0.9n]$, respectively. The results at the 0.1, 0.05 and 0.01 significance levels are reported in Table 1. When $n = 250$, the size is very close to the nominal 0.01 level and is acceptable at the nominal 0.05 level, but is quite conservative at the nominal 0.1 level. When $n$ is increased to 400, all size values are close to the nominal levels. Power increases when $n$ increases from 250 to 400. When $k = [0.9n]$, the power is lower than when $k = [0.5n]$. We also have the simulation results when $n = 200$. But in this case, all size values are small and power is very low, and hence, they are not reported here.

**5. Proofs of Theorems 2.1 and 2.2.** This section gives the proofs of Theorems 2.1 and 2.2.

PROOF OF THEOREM 2.1. Let $S_n = \sum_{t=1}^{n} X_t$ and $p = 1 + \iota$. From $K = 1, 2, \ldots$, let $l = [\sqrt{K}]$ and define $A_{t,K} = \|\sum_{j=1}^{K}[X_{t+j} - E(X_{t+j}|\mathcal{F}_l(t+j))]\|_p$ and $B_{t,K} = \|\sum_{j=1}^{K} E(X_{t+j}|\mathcal{F}_l(t+j))\|_p$. Since $E(X_{t+j}|\mathcal{F}_l(t+j))$ are $l$-dependent and the $L^p(\nu)$ NED assumption holds, it can be readily shown that there is a constant $\alpha > 0$ such that $A_{t,K} + B_{t,K} = O(K^{1-\alpha})$ uniformly in $t$. So, $\|\sum_{i=1}^{K} X_{t+i}\|_p = O(K^{(1-\alpha)})$. By Proposition 1 of [37], we have

$$\left\|\max_{1 \leq t \leq 2^k} |S_t|\right\|_p \leq \sum_{r=0}^{k} \left[\sum_{i=1}^{2^{k-r}} \|S_{2^r i} - S_{2^r(i-1)}\|_p^p\right]^{1/p}.$$

Thus, for some $0 < \rho < 1$, $\|\max_{1 \leq t \leq 2^k} |S_t|\|_p = O(2^{k\rho})$, from which (a) and (b) follow easily. □

TABLE 1
*Size and power of $\widehat{W}_n(1)$ for testing change-point in FARIMA$(0,d,0)$ models (1000 replications)*

|          | $n=250$ | | | $n=400$ | | |
|----------|--------|--------|--------|--------|--------|--------|
|          | **10%** | **5%** | **1%** | **10%** | **5%** | **1%** |
| $d_0$    | | Sizes | | | | |
| 0.1      | 0.055 | 0.039 | 0.012 | 0.081 | 0.049 | 0.015 |
| 0.2      | 0.059 | 0.037 | 0.012 | 0.083 | 0.046 | 0.014 |
| 0.3      | 0.064 | 0.038 | 0.010 | 0.078 | 0.047 | 0.012 |
| 0.4      | 0.050 | 0.031 | 0.010 | 0.077 | 0.041 | 0.014 |
| $d_{10}$ | Power when $d_0=0.1$ and $k=[0.5n]$ | | | | | |
| 0.2      | 0.168 | 0.120 | 0.040 | 0.304 | 0.235 | 0.126 |
| 0.3      | 0.333 | 0.260 | 0.114 | 0.655 | 0.566 | 0.403 |
| 0.4      | 0.658 | 0.571 | 0.387 | 0.924 | 0.886 | 0.791 |
| $d_{10}$ | Power when $d_0=0.1$ and $k=[0.9n]$ | | | | | |
| 0.2      | 0.135 | 0.089 | 0.022 | 0.180 | 0.114 | 0.056 |
| 0.3      | 0.181 | 0.124 | 0.056 | 0.303 | 0.225 | 0.106 |
| 0.4      | 0.424 | 0.334 | 0.197 | 0.582 | 0.498 | 0.312 |

To prove Theorem 2.2, we need the following lemma which is used for (5.2), (5.7) and (5.10).

LEMMA 5.1. *Let $X_t$ be defined as in Theorem 2.2. Then* (a)

$$\left\|\sum_{t=i+1}^{j}\{X_t - E[X_t|\mathcal{F}_{t-i}(t)]\}\right\|_{2+\iota} \Big/ \left[\sum_{t=i+1}^{j}\frac{1}{(t-i)^{2\nu}}\right]^{1/2} = O(1),$$

*uniformly in $j$ and $i < j$, and* (b) *furthermore, if (2.1) holds, then we have*

$$\left\|\sum_{t=1}^{j-1}\{E[X_{-t}|\mathcal{F}_{-t+j}(-t)] - E[X_{-t}|\mathcal{F}_{-t+j-1}(-t)]\}\right\|_{2+\iota_1} = O(1).$$

PROOF. Let $p = 2 + \iota$ and $\xi_{t,k} = X_t - E[X_t|\mathcal{F}_k(t)]$ for $k \geq 0$. Since $X_t$ is an $\mathcal{F}_t$-measurable martingale difference and $\{\varepsilon_t\}$ is independent, we know that $\xi_{t,t-i}$ is an $\mathcal{F}_t$-measurable martingale difference. By Definition 2.1, $\sup_{i<j}\sup_{i<t\leq j}[\|\xi_{t,t-i}\|_p(t-i)^\nu] = \sup_{i<j}\sup_{i<t\leq j}\sup_{0\leq k<\infty}(\|\xi_{t,k}\|_p k^\nu) \leq \sup_{0\leq k<\infty}\sup_{-\infty<t<\infty}(\|\xi_{t,k}\|_p k^\nu) = O(1)$. By Burkholder's inequality in [10], page 384, there exists a constant $B$, depending only on $\iota$, such that

$$E\left|\sum_{t=i+1}^{j}\xi_{t,t-i}\right|^p \leq BE\left(\sum_{t=i+1}^{j}|\xi_{t,t-i}|^2\right)^{p/2} \leq B\left(\sum_{t=i+1}^{j}\|\xi_{t,t-i}\|_p^2\right)^{p/2},$$

where the last step uses Minkowski's inequality. Thus, (a) holds. Since $E[X_{-t}|\ \mathcal{F}_{-t+j}(-t)] - E[X_{-t}|\mathcal{F}_{-t+j-1}(-t)]$ is an $\mathcal{F}_{-t}$-measurable martingale difference and $2\nu_1 > 1$, similarly, we can prove that (b) holds. $\square$

PROOF OF THEOREM 2.2. By Theorem 2 in [13], the proof of (a) is much easier than that of (b). So, only the latter is presented here.

Let $X_{t,i}^{(0)} = E[X_t|\mathcal{F}_{-i+1}(0)] - E[X_t|\mathcal{F}_{-i}(0)]$ for $i \leq -1$. Note that $E(X_t|\mathcal{F}_{-t}(0)) = EX_t = 0$ when $t \leq -1$. We have the decomposition

$$
\begin{aligned}
\sum_{t=-k}^{-1} X_t &= \sum_{t=-k}^{-1} \{X_t - E[X_t|\mathcal{F}_{k+1}(0)]\} + \sum_{t=-k}^{-1} \sum_{i=-k}^{t} X_{t,i}^{(0)} \\
&= \sum_{t=-k}^{-1} \{X_t - E[X_t|\mathcal{F}_{k+1}(0)]\} + \sum_{i=-k}^{-1} \sum_{t=i}^{-1} X_{t,i}^{(0)} \\
&= \sum_{t=-k}^{-1} \{X_t - E[X_t|\mathcal{F}_{k+1}(0)]\} + \sum_{j=1}^{k} \sum_{t=1}^{j} X_{-t,-j}^{(0)}.
\end{aligned}
\tag{5.1}
$$

Note that $E[X_t|\mathcal{F}_{k+1}(0)] = E[X_t|\mathcal{F}_{t+k+1}(t)]$ when $t \leq 0$ and $t+k+1 \geq 0$. By Lemma 5.1(a), $\|\sum_{t=-k}^{-1}\{X_t - E[X_t|\mathcal{F}_{k+1}(0)]\}\|_{2+\iota} = O([\sum_{t=-k}^{-1}(t+k+1)^{-2\nu}]^{1/2}) = O[(\sum_{t=1}^{k} t^{-2\nu})^{1/2}]$. Thus, by the Cauchy–Schwarz inequality, for any $\epsilon > 0$, we have

$$
\begin{aligned}
&P\left(\max_{l \leq k}\left|\frac{1}{k^{1/2-\delta}}\sum_{t=-k}^{-1}\{X_t - E[X_t|\mathcal{F}_{k+1}(0)]\}\right| > \epsilon\right) \\
&\leq \sum_{k=l}^{\infty} P\left(\frac{1}{k^{1/2-\delta}}\left|\sum_{t=-k}^{-1}\{X_t - E[X_t|\mathcal{F}_{k+1}(0)]\}\right| > \epsilon\right) \\
&\leq \frac{1}{\epsilon^{2+\iota}} \sum_{k=l}^{\infty} \frac{1}{k^{(1/2-\delta)(2+\iota)}} E\left|\sum_{t=-k}^{-1}\{X_t - E[X_t|\mathcal{F}_{k+1}(0)]\}\right|^{2+\iota} \\
&\leq \frac{O(1)}{\epsilon^{2+\iota}} \sum_{k=l}^{\infty} \frac{1}{k^{(1/2-2\delta)(2+\iota)}} \left(\sum_{t=1}^{k} \frac{1}{t^{2\nu+2\delta}}\right)^{1+\iota/2} \\
&= O\left(\frac{1}{l^{(1/2-2\delta)(2+\iota)-1}}\right),
\end{aligned}
\tag{5.2}
$$

as $\delta > 0$ is small enough such that $(1/2 - 2\delta)(2+\iota) > 1$. By Lemma 1 in [10], page 31 and (5.2),

$$
\sum_{t=-k}^{-1} \{X_t - E[X_t|\mathcal{F}_{k+1}(0)]\} = O(k^{1/2-\delta}) \qquad \text{a.s.}
\tag{5.3}
$$

The second term in (5.1) can be rewritten as

$$\sum_{j=2}^{k+1} Y_{0j} \quad \text{and} \quad Y_{0j} = \sum_{t=1}^{j-1} X^{(0)}_{-t,-j+1}.$$

By (5.3), it is sufficient for (b) to show that we can define a sequence of i.i.d. $m \times 1$ normal vectors $\{G_{2t}\}$ with mean zero and covariance $\Omega$ such that

(5.4) $$\sum_{j=2}^{k+1} Y_{0j} - \sum_{j=1}^{k} G_{2j} = O(k^{1/2-\delta}) \quad \text{a.s.}$$

Since $X^{(0)}_{-t,-j+1} \in \mathcal{F}_j(0)$, we know that $Y_{0j} \in \mathcal{F}_j(0)$ and $E(Y_{0j}|\mathcal{F}_{j-1}(0)) = 0$. Thus, $\{Y_{0j}, \mathcal{F}_j(0), j = 1, 2, \ldots\}$ is a sequence of forward martingale differences. Using the strong invariance principle in [13], Theorem 1, it is sufficient for (5.4) to verify the following conditions:

(i) there exists an $\tilde{\iota} > 0$ such that $E|Y_{0j}|^{2+\tilde{\iota}} \le M$, a constant, uniformly in $j$;

(ii) for some $\delta > 0$, the following holds uniformly in $s$:

(5.5) $$\frac{1}{n^{1-\delta}} \sum_{j=s+1}^{s+n} [E(Y_{0j}Y'_{0j}) - \Omega] = O(1),$$

(5.6) $$\frac{1}{n^{1-\delta}} E \left| \sum_{j=s+1}^{s+n} [E(Y_{0j}Y'_{0j}|\mathcal{F}_s(0)) - E(Y_{0j}Y'_{0j})] \right| = O(1).$$

Note that $X^{(0)}_{-t,-j+1} = E[X_{-t}|\mathcal{F}_j(0)] - E[X_{-t}|\mathcal{F}_{j-1}(0)]$ and $E[X_{-t}|\mathcal{F}_j(0)] = E[X_{-t}|\mathcal{F}_{-t+j}(-t)]$ when $t \ge 1$ and $-t+j \ge 0$. When $2\nu > 1$, by Minkowski's inequality and Lemma 5.1(a), it follows that, uniformly in $j$,

$$E|Y_{0j}|^{2+\iota} = E\left|\sum_{t=1}^{j-1} X^{(0)}_{-t,-j+1}\right|^{2+\iota}$$

(5.7) $$\le O(1) E\left|\sum_{t=1}^{j-1} \{X_{-t} - E[X_{-t}|\mathcal{F}_j(0)]\}\right|^{2+\iota}$$

$$+ O(1) E\left|\sum_{t=1}^{j-1} \{X_{-t} - E[X_{-t}|\mathcal{F}_{j-1}(0)]\}\right|^{2+\iota} = O(1).$$

That is, (i) holds. When $2\nu = 1$ and (2.1) is satisfied, (i) holds by Lemma 5.1(b).

For (ii), we make a decomposition as

$$Y_{0j}Y'_{0j} = \left(\sum_{t=1}^{s-1} X^{(0)}_{-t,-j+1}\right)\left(\sum_{t=1}^{s-1} X^{(0)}_{-t,-j+1}\right)'$$

$$+ \left(\sum_{t=s}^{j-1} X^{(0)}_{-t,-j+1}\right)\left(\sum_{t=s}^{j-1} X^{(0)}_{-t,-j+1}\right)'$$

(5.8)
$$+ \left[\left(\sum_{t=s}^{j-1} X^{(0)}_{-t,-j+1}\right)\left(\sum_{t=1}^{s-1} X^{(0)}_{-t,-j+1}\right)'\right.$$

$$\left.+ \left(\sum_{t=1}^{s-1} X^{(0)}_{-t,-j+1}\right)\left(\sum_{t=s}^{j-1} X^{(0)}_{-t,-j+1}\right)'\right]$$

$$\equiv B_{1js} + B_{2js} + B_{3js}.$$

We first show that

(5.9) $$\frac{1}{n^{1-\delta}} \sum_{j=s+1}^{s+n} (EB_{2js} - \Omega) = o(1),$$

uniformly in $s$. Let $Z_{t,j} = E[X_{-t}|\mathcal{F}_j(0)]$. Since $E(X^{(0)}_{-t,-j+1}|\mathcal{F}_{-t-1}) = 0$,

$$\sum_{j=s+1}^{s+n} EB_{2js} = \sum_{j=s+1}^{s+n} \sum_{t=s}^{j-1} E(X^{(0)}_{-t,-j+1} X^{(0)'}_{-t,-j+1})$$

$$= \sum_{j=s+1}^{s+n} \sum_{t=s}^{j-1} E(Z_{t,j} Z'_{t,j}) - \sum_{j=s+1}^{s+n} \sum_{t=s}^{j-1} E(Z_{t,j-1} Z'_{t,j-1})$$

$$= \sum_{j=s+1}^{s+n} \sum_{t=s}^{j-1} E(Z_{t,j} Z'_{t,j}) - \sum_{j=s}^{s+n-1} \sum_{t=s}^{j} E(Z_{t,j} Z'_{t,j})$$

$$= \sum_{t=s}^{s+n-1} E(Z_{t,s+n} Z'_{t,s+n}) - E(Z_{s,s} Z'_{s,s}) - \sum_{j=s+1}^{n+s-1} E(Z_{j,j} Z'_{j,j})$$

$$= \sum_{t=s}^{s+n-1} E(Z_{t,s+n} Z'_{t,s+n}),$$

where $Z_{j,j} = E(X_{-j}|\mathcal{F}_j(0)) = 0$ is used since $X_{-j}$ is independent of $\mathcal{F}_j(0)$. Since $Z_{t,j} = E[X_{-t}|\mathcal{F}_{-t+j}(-t)]$, by the near-epoch dependence of $X_t$,

$$\frac{1}{n^{1-\delta}} \left|\sum_{t=s}^{s+n-1} E(Z_{t,s+n} Z'_{t,s+n} - X_{-t} X'_{-t})\right|$$

$$\leq \frac{1}{n^{1-\delta}} \left[\sum_{t=s}^{s+n-1} E|X_{-t} - Z_{t,s+n}|^2 + 2\sum_{t=s}^{s+n-1} E(|X_{-t} - Z_{t,s+n}||X_{-t}|)\right]$$

$$\leq O\left(\frac{1}{n^{1-\delta}}\right)\left[\sum_{t=s}^{s+n-1} \frac{1}{(n+s-t)^{2\nu}} + \sum_{t=s}^{s+n-1} (E|X_{-t} - Z_{t,s+n}|^2 E|X_{-t}|^2)^{1/2}\right]$$

$$= O\left(\frac{1}{n^{1-\delta}}\right) \sum_{t=s}^{s+n-1} \frac{1}{(n+s-t)^\nu} = o(1),$$

for $\delta < \nu$, where $O(1)$ and $o(1)$ hold uniformly in $s$. Thus, (5.9) holds.

Since $X_{-t,-j+1}^{(0)}$ is an $\mathcal{F}_{-t}$-measurable martingale difference, we have $E|B_{1js}| \leq E|\sum_{t=1}^{s-1} X_{-t,-j+1}^{(0)}|^2 = \sum_{t=1}^{s-1} E|X_{-t,-j+1}^{(0)}|^2$. When $2\nu > 1$, by the near-epoch dependence of $X_t$, it follows that, uniformly in $j \geq s$,

$$\sum_{t=1}^{s-1} E|X_{-t,-j+1}^{(0)}|^2 \leq O(1) \sum_{t=1}^{s-1} \{E|X_{-t} - Z_{t,j}|^2 + E|X_{-t} - Z_{t,j-1}|^2\}$$

$$\leq O(1) \sum_{t=1}^{s-1} \frac{1}{(j-t)^{2\nu}}.$$

When $2\nu = 1$, by (2.1) and Minkowski's inequality, $E|X_{-t,-j+1}^{(0)}|^2 \leq O(1)(j-t)^{-2\nu_1}$ and $2\nu_1 > 1$. Letting $\tilde{\nu} = 2\nu$ or $2\nu_1$ according as $2\nu > 1$ or $2\nu = 1$, we have, for $0 < \delta < (\tilde{\nu}-1)/2$,

$$\frac{1}{n^{1-\delta}} \sum_{j=s+1}^{n+s} E|B_{1js}| \leq \frac{O(1)}{n^{1-\delta}} \sum_{j=s+1}^{s+n} \sum_{t=1}^{s-1} \frac{1}{(j-t)^{\tilde{\nu}}}$$

(5.10) $$\leq O(1) \frac{1}{n^\delta} \sum_{j=s+1}^{s+n} \frac{1}{(j-s)^{1-2\delta+\tilde{\delta}}} \sum_{t=1}^{s-1} \frac{1}{(s-t)^{\tilde{\nu}-\tilde{\delta}}}$$

$$= O(1),$$

uniformly in $s$ as $n \to \infty$, for $2\delta < \tilde{\delta} < \tilde{\nu} - 1$, where we have used $j - s \leq \min\{j-t, n\}$ and $j - t \geq s - t$. By (5.9)–(5.10) and the Cauchy–Schwarz inequality, we can show that

(5.11) $$\frac{1}{n^{1-\delta}} \sum_{j=s+1}^{n+s} E|B_{3js}| = O(1).$$

By (5.8)–(5.11), we can establish (5.5). Since $B_{2js} \in \mathcal{F}_{-s}$ is independent of $\mathcal{F}_s(0)$ when $j > s$, $E[B_{2js}|\mathcal{F}_s(0)] = EB_{2js}$. By (5.8)–(5.11), uniformly in $s$, we have

(5.12) $$\frac{1}{n^{1-\delta}} E \left| \sum_{j=s+1}^{n+s} \{E[Y_{0j}Y_{0j}'|\mathcal{F}_s(0)] - \Omega\} \right| = O(1).$$

By (5.5) and (5.12), (5.6) holds. □

**6. Proof of Theorem 3.1.** We first present three lemmas. Lemma 6.1 comes directly from Theorem 2.2, while Lemma 6.2 can be proved by using Theorem 2.1 and the details are given in [28].

LEMMA 6.1. *If Assumptions* 3.1(ii) *and* 3.2(ii) *hold, then in the sense of Theorem* 2.2, *we can define i.i.d.* $m \times 1$ *normal vector sequences,* $\{G_{1t}\}$ *and* $\{G_{2t}\}$, *with mean zero and covariance* $\Omega$ *such that, for some* $\delta > 0$,

$$\text{(a)} \quad \max_{g_n \le k \le n} k^\delta \left| \frac{1}{\sqrt{k}} \sum_{t=1}^{k} D_t(\lambda_0) - \frac{1}{\sqrt{k}} \sum_{t=1}^{k} G_{1t} \right| = o_p(1),$$

$$\text{(b)} \quad \max_{g_n \le k < n} k^\delta \left| \frac{1}{\sqrt{k}} \sum_{t=-k}^{-1} D_t(\lambda_0) - \frac{1}{\sqrt{k}} \sum_{t=-k}^{-1} G_{2t} \right| = o_p(1).$$

LEMMA 6.2. *Let* $\Theta_0(k, \eta) = \{\lambda : k^{\tilde{\iota}} |\lambda - \lambda_0| < \eta\}$. *Suppose Assumptions* 3.1(iii)–(iv), 3.2(iii) *and* 3.3(iii) *hold. For any* $\epsilon > 0$, (1) *if* $\tilde{\iota} = 0$, *then there is* $\eta > 0$ *such that the following* (a)–(b) *hold with* $\delta = 0$, *and* (2) *if* $\tilde{\iota} > 0$, *then there is* $\delta > 0$ *such that the following* (a)–(b) *hold for any fixed* $\eta > 0$:

$$\text{(a)} \quad \lim_{n \to \infty} P\left( \max_{g_n \le k \le n} \max_{\Theta_0(k,\eta)} \frac{1}{k^{1-\delta}} \left| \sum_{t=1}^{k} [\widetilde{P}_t(\lambda) - \Sigma] \right| \ge \epsilon \right) = 0,$$

$$\text{(b)} \quad \lim_{n \to \infty} P\left( \max_{g_n \le n-k < n} \max_{\Theta_0(n-k,\eta)} \frac{1}{(n-k)^{1-\delta}} \left| \sum_{t=k+1}^{n} [\widetilde{P}_t(\lambda) - \Sigma] \right| \ge \epsilon \right) = 0.$$

LEMMA 6.3. *If the assumptions of Theorem* 3.1 *hold, then there exists a* $\delta > 0$ *such that* $\hat{\lambda}_n(k)$ *and* $\hat{\lambda}_{1n}(k)$ *have the uniform expansions*

$$\text{(a)} \quad \max_{g_n \le k \le n} k^\delta \left| \sqrt{k}[\hat{\lambda}_n(k) - \lambda_0] - \frac{\Sigma^{-1}}{\sqrt{k}} \sum_{t=1}^{k} D_t(\lambda_0) \right| = o_p(1),$$

$$\text{(b)} \quad \max_{g_n \le n-k < n} (n-k)^\delta \left| \sqrt{n-k}[\hat{\lambda}_{1n}(k) - \lambda_0] - \frac{\Sigma^{-1}}{\sqrt{n-k}} \sum_{t=k+1}^{n} D_t(\lambda_0) \right| = o_p(1).$$

PROOF. We only prove part (b). By Lemma A.1(b) in the Appendix,

$$P\left( \max_{g_n \le n-k < n} |\hat{\lambda}_{1n}(k) - \lambda_0| > \epsilon \right)$$

$$= P\left\{ |\hat{\lambda}_{1n}(k) - \lambda_0| > \epsilon, \sum_{t=k+1}^{n} [\tilde{l}_t(\hat{\lambda}_{1n}(k)) - \tilde{l}_n(\lambda_0)] \ge 0, \right.$$

$$\left. \text{for some } k \in [1, n - g_n] \right\}$$

$$\leq P\left\{\max_{g_n\leq n-k<n}\sup_{|\lambda-\lambda_0|>\epsilon}\sum_{t=k+1}^n [\tilde{l}_t(\lambda)-\tilde{l}_t(\lambda_0)]\geq 0\right\}=o(1),$$

for any $\epsilon>0$ and as $n\to\infty$. Thus,

(6.1) $$\max_{g_n\leq n-k<n}|\hat{\lambda}_{1n}(k)-\lambda_0|=o_p(1).$$

Using Taylor's expansion for each element of $\partial L_{1n}(k,\hat{\lambda}_{1n}(k))/\partial\lambda=0$, we have

(6.2) $$\hat{\lambda}_{1n}(k)-\lambda_0=\left(\frac{1}{n-k}\sum_{t=k+1}^n \widetilde{P}^*_{nt}\right)^{-1}\frac{1}{n-k}\sum_{t=k+1}^n \widetilde{D}_t(\lambda_0),$$

for each $k$, where the $i$th row of $\widetilde{P}^*_{nt}$ is the $i$th row of $\widetilde{P}_t(\hat{\lambda}^{*(i)}_{1n}(k))$ for some $\hat{\lambda}^{*(i)}_{1n}(k)$ such that $|\hat{\lambda}^{*(i)}_{1n}(k)-\lambda_0|\leq|\hat{\lambda}_{1n}(k)-\lambda_0|$ for $i=1,\ldots,m$. Observing that $D_t(\lambda_0)$ is strictly stationary, by Lemma 6.1(b), the law of iterated logarithm (LIL) and Assumption 3.3(ii), it follows that, for any $\delta_0>0$,

(6.3) $$\max_{g_n\leq n-k<n}\left|\frac{1}{(n-k)^{0.5+\delta_0}}\sum_{t=k+1}^n \widetilde{D}_t(\lambda_0)\right|=o_p(1).$$

Let $\delta_1\in(0,1/2)$. By Lemma 6.2(b) with $\tilde{\iota}=0$ and (6.1)–(6.3),

(6.4) $$\max_{g_n\leq n-k<n}|(n-k)^{\delta_1}[\hat{\lambda}_{1n}(k)-\lambda_0]|=o_p(1).$$

By (6.2)–(6.4) and Lemma 6.2(b) with $\tilde{\iota}=\delta_1$, there exists some $\delta>0$ such that

$$\max_{g_n\leq n-k<n}(n-k)^\delta\left|\sqrt{n-k}[\hat{\lambda}_{1n}(k)-\lambda_0]-\frac{\Sigma^{-1}}{\sqrt{n-k}}\sum_{t=k+1}^n \widetilde{D}_t(\lambda_0)\right|$$

$$=\max_{g_n\leq n-k<n}\left|\left(\frac{1}{n-k}\sum_{t=k+1}^n \widetilde{P}^*_{nt}\right)^{-1}\left[\frac{1}{(n-k)^{1-2\delta}}\sum_{t=k+1}^n (\Sigma-\widetilde{P}^*_{nt})\right]\right.$$

$$\left.\times\frac{\Sigma^{-1}}{(n-k)^{1/2+\delta}}\sum_{t=k+1}^n \widetilde{D}_t(\lambda_0)\right|=o_p(1).$$

Furthermore, by Assumption 3.3(ii), (b) holds. □

We further need two lemmas, which are directly used for Theorem 3.1.

LEMMA 6.4. *Under the assumptions of Theorem* 3.1, *it follows that*

(a) $\max_{k\in\Pi_n}|W_n(k)-S_n(k)|=o_p(1)$  *and*  (b) $\max_{k\notin\Pi_n}W_n(k)=O_p(g_n)$,

where $\Pi_n = [\log n, n - \log n]$ and

$$S_n(k) = \left|\frac{\Omega^{-1/2}}{\sqrt{k}}\sum_{t=1}^{k} D_t(\lambda_0)\right|^2 + \left|\frac{\Omega^{-1/2}}{\sqrt{n-k}}\sum_{t=k+1}^{n} D_t(\lambda_0)\right|^2 - \left|\frac{\Omega^{-1/2}}{\sqrt{n}}\sum_{t=1}^{n} D_t(\lambda_0)\right|^2.$$

PROOF. (a) By Lemma 6.3, we have

(6.5) $$\max_{k\in\Pi_n}\left|\sqrt{k}[\hat{\lambda}_n(k) - \lambda_0] - \frac{\Sigma^{-1}}{\sqrt{k}}\sum_{t=1}^{k} D_t(\lambda_0)\right| = O_p(\log^{-\delta} n),$$

(6.6) $$\max_{k\in\Pi_n}\left|\sqrt{n-k}[\hat{\lambda}_{1n}(k) - \lambda_0] - \frac{\Sigma^{-1}}{\sqrt{n-k}}\sum_{t=k+1}^{n} D_t(\lambda_0)\right| = O_p(\log^{-\delta} n),$$

for some $\delta > 0$. As for (6.3), by Lemma 6.1 and the LIL, we can show that

(6.7) $$\max_{k\in\Pi_n}|\sqrt{k}[\hat{\lambda}_n(k) - \lambda_0]| = O_p[(\log\log n)^{1/2}],$$

(6.8) $$\max_{k\in\Pi_n}|\sqrt{n-k}[\hat{\lambda}_{1n}(k) - \lambda_0]| = O_p[(\log\log n)^{1/2}].$$

By Lemma 6.2(a)–(b) and Lemma 6.3, we know that $\max_{k\in\Pi_n}|\hat{\Sigma}_n(k)/n - \Sigma| = O_p(n^{-\delta})$. By Lemma A.2 in the Appendix, we have $\max_{k\in\Pi_n}|\hat{\Omega}_n(k)/n - \Omega| = O_p(n^{-\delta})$. Furthermore, by (6.7)–(6.8), it follows that

(6.9) $$\max_{k\in\Pi_n}\left|W_n(k) - \frac{k(n-k)}{n}[\hat{\lambda}_n(k) - \hat{\lambda}_{1n}(k)]'\Omega_0[\hat{\lambda}_n(k) - \hat{\lambda}_{1n}(k)]\right| = o_p(1),$$

where $\Omega_0 = \Sigma\Omega^{-1}\Sigma$. Denote

$$\xi_n(k) = \sqrt{\frac{k(n-k)}{n}}\left[\frac{1}{k}\sum_{t=1}^{k} D_t(\lambda_0) - \frac{1}{n-k}\sum_{t=k+1}^{n} D_t(\lambda_0)\right].$$

By (6.5)–(6.6), we have

(6.10) $$\max_{k\in\Pi_n}\left|\sqrt{\frac{k(n-k)}{n}}[\hat{\lambda}_n(k) - \hat{\lambda}_{1n}(k)] - \Sigma^{-1}\xi_n(k)\right| = O_p(\log^{-\delta} n).$$

By (6.7)–(6.8) and (6.10), it follows that

(6.11) $$\max_{k\in\Pi_n}\left|\left\{\sqrt{\frac{k(n-k)}{n}}[\hat{\lambda}_n(k) - \hat{\lambda}_{1n}(k)] - \Sigma^{-1}\xi_n(k)\right\}'\Omega_0 \times \left\{\sqrt{\frac{k(n-k)}{n}}[\hat{\lambda}_n(k) - \hat{\lambda}_{1n}(k)]\right\}\right| = o_p(1).$$

By (6.9) and (6.11), we can show that $\max_{k\in\Pi_n}|W_n(k) - \xi_n'(k)\Omega^{-1}\xi_n(k)| = o_p(1)$. By direct calculation, we have $S_n(k) = \xi_n'(k)\Omega^{-1}\xi_n(k)$. Thus, (a) holds. The proof of (b) is easy and can be found in [28]. □

LEMMA 6.5. *Let $\{G_t, t = 1, 2, \ldots\}$ be an i.i.d. sequence of $m \times 1$ random vectors with $EG_t = 0$ and $E(G_t G_t') = I$. If $E|G_t|^{2+\iota} < \infty$ for some $\iota > 0$, then, for each $\mu \in (0, 1)$ and for any $x$,*

$$P\left(\frac{1}{a_n(m)}\left[\max_{\log n \leq k \leq \mu n}\left|\frac{1}{\sqrt{k}}\sum_{t=1}^k G_t\right|^2 - b_n(m)\right] \leq x\right) \to \exp(-e^{-x/2}),$$

*as $n \to \infty$, where $a_n(m)$ and $b_n(m)$ are defined in Theorem 3.1.*

PROOF. The lemma can be proved readily by using Lemma 2.2 of [17] and Corollary A.2 of [12]. □

PROOF OF THEOREM 3.1. Let $S_n(k)$ be defined as in Lemma 6.4 and denote

$$\widetilde{S}_1(k) = \frac{\Omega^{-1/2}}{\sqrt{k}}\sum_{t=1}^k D_t(\lambda_0) \quad \text{and} \quad \widetilde{S}_{2n}(k) = \frac{\Omega^{-1/2}}{\sqrt{n-k}}\sum_{t=k+1}^n D_t(\lambda_0).$$

Let $\mu \in (0, 0.5)$. By Lemma 6.1(a) and the continuous mapping theorem,

$$\left|\max_{\log n \leq k \leq \mu n} S_n(k) - \max_{\log n \leq k \leq \mu n}|\widetilde{S}_1(k)|^2\right|$$
$$\leq \max_{\log n \leq k \leq \mu n}|S_n(k) - |\widetilde{S}_1(k)|^2|$$
$$= \max_{\log n \leq k \leq \mu n}||\widetilde{S}_{2n}(k)|^2 - |\widetilde{S}_1(n)|^2|$$
$$\to_{\mathcal{L}} \max_{0 \leq \tau \leq \mu}\left|\frac{|B(1) - B(\tau)|^2}{1 - \tau} - |B(1)|^2\right|,$$

as $n \to \infty$, where $\to_{\mathcal{L}}$ denotes convergence in distribution and $\{B(\tau): \tau \in [0, 1]\}$ is a standard Brownian motion. Thus, for any $\epsilon > 0$,

(6.12) $\quad \limsup_{n \to \infty} P\left(\left|\max_{\log n \leq k \leq \mu n} S_n(k) - \max_{\log n \leq k \leq \mu n}|\widetilde{S}_1(k)|^2\right| > \epsilon\right) \to 0,$

as $\mu \to 0$. Similarly, as $\mu \to 0$,

(6.13) $\quad \limsup_{n \to \infty} P\left(\left|\max_{\log n \leq n-k \leq \mu n} S_n(k) - \max_{\log n \leq n-k \leq \mu n}|\widetilde{S}_{2n}(k)|^2\right| > \epsilon\right) \to 0.$

Denote $B_1(k) = \Omega^{-1/2}\sum_{t=1}^k G_{1t}/\sqrt{k}$ and $B_2(k) = \Omega^{-1/2}\sum_{t=-k}^{-1} G_{2t}/\sqrt{k}$, where $\{G_{1t}\}$ and $\{G_{2t}\}$ are defined as in Lemma 6.1. By Lemma 6.1(a), for each $\mu$, we have

$$\left|\max_{\log n \leq k \leq \mu n}|\widetilde{S}_1(k)|^2 - \max_{\log n \leq k \leq \mu n}|B_1(k)|^2\right|$$
$$\leq \max_{\log n \leq k < n}|\sqrt{\log \log k}[\widetilde{S}_1(k) - B_1(k)]|\frac{|B_1(k)| + |\widetilde{S}_1(k)|}{\sqrt{\log \log k}} = o_p(1),$$

as $n \to \infty$. Furthermore, by Lemma 6.5(a), for each $\mu$ and $x$, it follows that

$$(6.14) \quad \lim_{n \to \infty} P\left(\left[\max_{\log n \leq k \leq \mu n} |\widetilde{S}_1(k)|^2 - b_n(m)\right] \Big/ a_n(m) < x\right) = \exp(-e^{-x/2}).$$

Applying the same argument to $\widetilde{S}_2^*(k) = (k\Omega)^{-1/2} \sum_{t=-k}^{1} D_t(\lambda_0)$ and $B_2(k)$ with the help of Lemma 6.1(b), and observing that $\max_{\log n \leq k \leq \mu n} |\widetilde{S}_2^*(k)|^2$ has the same distribution as $\max_{\log n \leq n-k \leq \mu n} |\widetilde{S}_{2n}(k)|^2$, we have

$$(6.15) \quad \lim_{n \to \infty} P\left(\left[\max_{\log n \leq n-k \leq \mu n} |\widetilde{S}_{2n}(k)|^2 - b_n(m)\right] \Big/ a_n(m) < x\right) = \exp(-e^{-x/2}).$$

Using a similar method as for (5.2), we can show that

$$\Delta_n \equiv \max_{\log n \leq n-k \leq \mu n} \left| \frac{\Omega^{-1/2}}{(n-k)^{0.5-\delta}} \sum_{t=k+1}^{n} \{D_t(\lambda_0) - E[D_t(\lambda_0)|\mathcal{F}_{t-k}(t)]\} \right| = o_p(1),$$

for some $\delta > 0$. Let $S_{2n}^0(k) = \Omega^{-1/2} \sum_{t=k+1}^{n} E[D_t(\lambda_0)|\mathcal{F}_{t-k}(t)]/\sqrt{n-k}$. So,

$$\left| \max_{\log n \leq n-k \leq \mu n} |\widetilde{S}_{2n}(k)|^2 - \max_{\log n \leq n-k \leq \mu n} |S_{2n}^0(k)|^2 \right|$$

$$\leq \max_{\log n \leq n-k \leq \mu n} ||\widetilde{S}_{2n}(k)|^2 - |S_{2n}^0(k)|^2|$$

$$\leq \Delta_n^2 + 2\Delta_n \max_{\log n \leq n-k \leq \mu n} \left| \frac{1}{(n-k)^{0.5+\delta}} \sum_{t=k+1}^{n} D_t(\lambda_0) \right| = o_p(1),$$

where the last step holds by Lemma 6.1(b), the LIL and the strict stationarity of $\{D_t(\lambda_0)\}$. Furthermore, by (6.15), for each $\mu$ and $x$, it follows that

$$\lim_{n \to \infty} P\left(\left[\max_{\log n \leq n-k \leq \mu n} |S_{2n}^0(k)|^2 - b_n(m)\right] \Big/ a_n(m) < x\right) = \exp(-e^{-x/2}).$$

By (6.14), the above two equations and independence of $\max_{\log n \leq n-k \leq \mu n} |S_{2n}^0(k)|^2$ and $\max_{\log n \leq k \leq \mu n} |\widetilde{S}_1(k)|^2$, for each $\mu \in (0, 0.5)$ and $x$, it follows that

$$P\left(\left[\max\left\{\max_{\log n \leq k \leq \mu n} |\widetilde{S}_1(k)|^2, \max_{\log n \leq n-k \leq \mu n} |\widetilde{S}_{2n}(k)|^2\right\} - b_n(m)\right] \Big/ a_n(m) < x\right)$$

$$= \exp(-2e^{-x/2}) + o(1).$$

Since $a_n(m) = 1 + o(1)$, by (6.12)–(6.13) and the preceding equation, we can show that, for each $x$ and any $\epsilon > 0$, there exist $N > 0$ and a constant

$\mu_0 \in (0, 1/2)$ such that, as $n > N$,

$$\left| P\left( \left[ \max\left\{ \max_{\log n \leq k \leq \mu_0 n} S_n(k), \max_{\log n \leq n-k \leq \mu_0 n} S_n(k) \right\} - b_n(m) \right] \middle/ a_n(m) < x \right) - \exp(-2e^{-x/2}) \right| < \frac{\epsilon}{2}.$$

By Lemma 6.1(a) and the continuous mapping theorem, we have

$$\max_{\mu_0 n \leq k \leq n-\mu_0 n} S_n(k) \longrightarrow_{\mathcal{L}} \max_{\mu_0 \leq \tau \leq 1-\mu_0} \left\{ \frac{|B(\tau)|^2}{\tau} + \frac{|B(1) - B(\tau)|^2}{1-\tau} - |B(1)|^2 \right\},$$

as $n \to \infty$. By the preceding two equations, for any $x$, we can show that

$$\lim_{n \to \infty} P\left( \left[ \max_{k \in [\log n, n - \log n]} S_n(k) - b_n(m) \right] \middle/ a_n(m) < x \right) = \exp(-2e^{-x/2}).$$

Finally, by Lemma 6.4(a)–(b), the conclusion holds. □

**7. Proof of Theorem 4.1.** It is sufficient for Theorem 4.1 to verify Assumptions 3.1–3.3. For simplicity, we only consider the case with $p = q = 0$, while the general case can be similarly verified.

In this case, the following expansions hold:

$$(7.1) \quad y_t = \sum_{i=0}^{\infty} c_{0i} \varepsilon_{t-i} \quad \text{and} \quad \varepsilon_t(\lambda) = (1-B)^d y_t = \sum_{i=0}^{\infty} a_{0i}(\lambda) y_{t-i},$$

where $c_{00} = a_{00}(\lambda) = 1$, $c_{0i} = O(i^{-1+d_0})$ and $a_{0i}(\lambda) = O(i^{-1-d})$. We further have

$$\frac{\partial \varepsilon_t(\lambda)}{\partial d} = \log(1-B)(1-B)^d y_t,$$

$$\frac{\partial^2 \varepsilon_t(\lambda)}{\partial d^2} = \log^2(1-B)(1-B)^d y_t \quad \text{and} \quad \frac{\partial^3 \varepsilon_t(\lambda)}{\partial d^3} = \log^3(1-B)(1-B)^d y_t,$$

where $\log(1-B) = -\sum_{i=1}^{\infty} B^i/i$ and $\log^k(1-B)(1-B)^d = \sum_{i=1}^{\infty} a_{ki} B^i$ with $a_{ki}(\lambda) = O(i^{-1-d} \log^k i)$ for $k = 1, 2, 3$.

PROOF OF ASSUMPTIONS 3.1–3.2. By Assumption 4.1, $\{y_t\}$ is strictly stationary with $E|y_t|^{2+\iota} < \infty$. Since $\Theta$ is compact, there exist constants $\underline{d}$ and $\tilde{d}$ such that $0 < \underline{d} \leq d \leq \tilde{d} < 0.5$ and $d_0 \in (\underline{d}, \tilde{d})$. Thus, $\sup_{\lambda \in \Theta} |a_{ki}(\lambda)| = O(i^{-1-\underline{d}} \log^k i)$ for $k = 0, 1, 2, 3$, and

$$\sup_{\Theta} |\varepsilon_t(\lambda)| = \sup_{\Theta} \left| \sum_{i=0}^{\infty} a_{0i}(\lambda) y_{t-i} \right| \leq |y_t| + O(1) \sum_{i=1}^{\infty} \frac{1}{i^{1+\underline{d}}} |y_{t-i}|.$$

Treating $\sup_\Theta |\varepsilon_t(\lambda)|$ and $y_t$ as elements in the $L^{2+\iota}$ space, we have

$$(7.2) \quad \left\|\sup_\Theta |\varepsilon_t(\lambda)|\right\|_{2+\iota} \leq O(1)\left[\|y_t\|_{2+\iota} + \sum_{i=1}^\infty \frac{\|y_{t-i}\|_{2+\iota}}{i^{1+\underline{d}}}\right] < \infty,$$

that is, the first part of Assumption 3.1(i) holds. The proof of the second part and Assumption 3.1(ii) can be found in [29]. Similar to (7.2), it can be proved that $\|\sup_\Theta[|\partial^2 \varepsilon_t(\lambda)/\partial d^2| + |\partial^3 \varepsilon_t(\lambda)/\partial d^3|]\|_{2+\iota} < \infty$. Thus, we can show that Assumption 3.1(iii)–(iv) holds. For a (large) integer $K$,

$$\left\|\sup_\Theta \left|\varepsilon_t(\lambda) - \sum_{i=0}^K a_{0i}(\lambda)y_{t-i}\right|\right\|_{2+\iota} \leq O(1) \sum_{i=K+1}^\infty \frac{\|y_{t-i}\|_{2+\iota}}{i^{1+\underline{d}}} = O\left(\frac{1}{K^{\underline{d}}}\right).$$

When $k > i$, by Lemma 2 in [35], it follows that

$$\|y_{t-i} - E[y_{t-i}|\mathcal{F}_k(t)]\|_{2+\iota}^{2+\iota} = E\left|\sum_{j=k-i}^\infty c_{0j}\varepsilon_{t-i-j}\right|^{2+\iota}$$

$$= O(1)\left(\sum_{j=k-i}^\infty c_{0j}^2\right)^{1+\iota/2} = O\left[\frac{1}{(k-i)^{(1-2d_0)(1+\iota/2)}}\right].$$

Let $K = [k/2]$. By the expansion of $\varepsilon_t(\lambda)$ in (7.1) and the preceding two equations,

$$\left\|\sup_\Theta |\varepsilon_t(\lambda) - E[\varepsilon_t(\lambda)|\mathcal{F}_k(t)]|\right\|_{2+\iota}^{2+\iota}$$

$$\leq O\left(\frac{1}{k^{\underline{d}(2+\iota)}}\right) + O(1)\left\{\sum_{i=0}^K \sup_\Theta |a_{0i}(\lambda)|\|y_{t-i} - E[y_{t-i}|\mathcal{F}_k(t)]\|_{2+\iota}\right\}^{2+\iota}$$

$$\leq O\left(\frac{1}{k^{\underline{d}(2+\iota)}}\right) + O(1)\left[\sum_{i=0}^K \frac{1}{(i+1)^{1+\underline{d}}(k-i)^{(1-2d_0)/2}}\right]^{2+\iota} = O\left(\frac{1}{k^{2\nu_0}}\right),$$

for some $\nu_0 > 0$. Using this with (7.2), we can show that Assumption 3.2(i) holds. Similarly, we can show that Assumption 3.2(iii) holds. Note that $(1-B)^{d_0}y_t = \varepsilon_t$. By Lemma 2 in [35],

$$\|D_t(\lambda_0) - E[D_t(\lambda_0)|\mathcal{F}_k(t)]\|_{2+\iota}^{2+\iota} = \sigma_\iota E\left|\sum_{i=k}^\infty \frac{\varepsilon_{t-i}}{i}\right|^{2+\iota} = O(k^{-(2+\iota)\times 0.5}),$$

where $\sigma_\iota = E|\varepsilon_t|^{2+\iota} < \infty$. Thus, $2\nu = 1$. Uniformly in $t$, it follows that

$$\|E[D_t(\lambda_0)|\mathcal{F}_{k+1}(t)] - E[D_t(\lambda_0)|\mathcal{F}_k(t)]\|_{2+\iota}^{2+\iota} = \sigma_\iota E\left|\frac{\varepsilon_{t-k}}{k}\right|^{2+\iota} = O\left(\frac{1}{k^{2+\iota}}\right).$$

Thus, we have that $\iota_1 = \iota$ and $\nu_1 = 1$. By the preceding two equations, we know that Assumption 3.2(ii) holds. □

PROOF OF ASSUMPTION 3.3. Since $y_t = 0$ for $t \leq 0$, by (7.1), we have

$$E\left[\sup_\Theta |\varepsilon_t(\lambda) - \tilde{\varepsilon}_t(\lambda)|\right]^2 = E\left[\sup_\Theta \left|\sum_{i=t}^\infty a_{0i}(\lambda) y_{t-i}\right|\right]^2 = O(t^{-2\underline{d}}).$$

Using this with (7.2), we can show that Assumption 3.3(i) holds. Similarly, we can show that the first part of Assumption 3.3(iii) holds. We now verify the second parts of Assumption 3.3(ii) and 3.3(iii). Denote

$$A_t = \varepsilon_t(\lambda_0) - \tilde{\varepsilon}_t(\lambda_0) = \sum_{i=t}^\infty a_{0i}(\lambda_0) y_{t-i},$$

$$A_{1t} = \frac{\partial \varepsilon_t(\lambda_0)}{\partial d} - \frac{\partial \tilde{\varepsilon}_t(\lambda_0)}{\partial d} = \sum_{i=t}^\infty a_{1i}(\lambda_0) y_{t-i},$$

$$A_{2t} = \frac{\partial \varepsilon_t(\lambda_0)}{\partial d} - v_t = -\sum_{i=t}^\infty \frac{1}{i} \varepsilon_{t-i},$$

where $v_t = -\sum_{i=1}^{t-1} \varepsilon_{t-i}/i$. We next make the decomposition

$$
(7.3) \quad
\begin{aligned}
-D_t(\lambda_0) + \tilde{D}_t(\lambda_0) &= \varepsilon_t(\lambda_0) A_{1t} + \frac{\partial \tilde{\varepsilon}_t(\lambda_0)}{\partial d} A_t \\
&= A_t A_{2t} - A_t A_{1t} + A_t v_t + \varepsilon_t(\lambda_0) A_{1t}.
\end{aligned}
$$

By (7.1), we can write $A_t$ as

$$A_t = \sum_{i=t}^\infty \sum_{j=i}^\infty a_{0i}(\lambda_0) c_{0j-i} \varepsilon_{t-j} = \sum_{j=t}^\infty \left[\sum_{i=t}^j a_{0i}(\lambda_0) c_{0j-i}\right] \varepsilon_{t-j}.$$

By Lemma 2 in [35], we can show that $E|A_t|^{2+\iota}$ is bounded by

$$C\left\{\sum_{j=t}^\infty \left[\sum_{i=t}^j a_{0i}(\lambda_0) c_{0j-i}\right]^2\right\}^{1+\iota/2}$$

$$\leq O(1)\left\{\sum_{j=t}^\infty \left[\sum_{i=t}^j \frac{1}{(j-i+1)^{1-d_0} i^{1+d_0}}\right]^2\right\}^{1+\iota/2},$$

where $C$ is some constant independent of $t$. Furthermore, by Lemma A.3 with $u = 1 - d_0$ and $v = 1 + d_0$, it follows that

$$(7.4) \quad E|A_t|^{2+\iota} \leq O(1)\left[\sum_{j=t}^\infty \left(\frac{1}{j^{1-d_0} t^{d_0}}\right)^2\right]^{1+\iota/2} = O(t^{-(1+\iota/2)}).$$

Furthermore, we can show that

(7.5)  $E|A_{1t}|^{2+\iota} = O((\log^2 t/t)^{1+\iota/2})$  and  $E|A_{2t}|^{2+\iota} = O(t^{-(1+\iota/2)})$.

By (7.3)–(7.5), we can show that $\|D_t(\lambda_0) - \widetilde{D}_t(\lambda_0)\|_{1+\iota/2} = O(t^{-1/2}\log t)$, that is, the second part of Assumption 3.3(iii) holds.

Denote $\tilde{k} = n - k$. By (7.4)–(7.5), $E|A_t A_{it}| \le (EA_t^2 EA_{it}^2)^{1/2} \le (E|A_t|^{2+\iota})^{1/(2+\iota)}(E|A_{it}|^{2+\iota})^{1/(2+\iota)} = O(t^{-1}\log^2 t)$. When $i = 1, 2$, by Lemma A.3 with $u = 1/2 - \delta$ and $v = 1$ and the Cauchy–Schwarz inequality, we can show that

(7.6)
$$P\left(\max_{g_n \le \tilde{k} < n} \frac{1}{\tilde{k}^{1/2-\delta}} \sum_{t=k+1}^{n} |A_t A_{it}| > \epsilon\right)$$
$$\le P\left(\sum_{t=1}^{n} \frac{|A_t A_{it}|}{(n-t+1)^{1/2-\delta}} > \epsilon\right) \to 0.$$

Next, consider the third term in (7.3). We first make the decomposition

$$\sum_{t=k+1}^{n} v_t A_t = -\sum_{i=1}^{k-1} \frac{1}{i} \sum_{t=k+1}^{n} \varepsilon_{t-i} A_t - \sum_{i=k}^{n-1} \frac{1}{i} \sum_{t=i+1}^{n} \varepsilon_{t-i} A_t,$$

where the last term is obtained by exchanging order of $\sum_{t=k+1}^{n} \sum_{i=k}^{t-1} \varepsilon_{t-i}A_t/i$.

Since $A_t$ is $\mathcal{F}_0$-measurable, by Lemma 2 in [35] and Minkowski's inequality, there exists a constant $B$, depending on $\iota$ and $E|\varepsilon_t|^{2+\iota}$, such that

$$\left\|\sum_{t=k+1}^{n} \varepsilon_{t-i} A_t\right\|_{2+\iota}^{2+\iota} \le B\left[\sum_{t=k+1}^{n}(E|A_t|^{2+\iota})^{2/(2+\iota)}\right]^{1+\iota/2} \le O(1)\left(\sum_{t=k+1}^{n} \frac{1}{t}\right)^{1+\iota/2},$$

uniformly in $i = 1, \ldots, n$. Let $\delta$ be small enough such that $(1/2 - 2\delta)(2+\iota) > 1$. Since $\tilde{k} \ge n - t + 1$, by the Markov and Minkowski inequalities, we have

$$P\left(\max_{g_n \le \tilde{k} < n} \frac{1}{\tilde{k}^{1/2-\delta}} \left|\sum_{i=1}^{k-1} \frac{1}{i} \sum_{t=k+1}^{n} \varepsilon_{t-i} A_t\right| > \epsilon\right)$$
$$\le O(1) \sum_{\tilde{k}=g_n}^{n-1} \frac{1}{\tilde{k}^{(1/2-\delta)(2+\iota)}} \left\|\sum_{i=1}^{k-1} \frac{1}{i} \sum_{t=k+1}^{n} \varepsilon_{t-i} A_t\right\|_{2+\iota}^{2+\iota}$$
$$\le O(1) \sum_{\tilde{k}=g_n}^{n-1} \frac{1}{\tilde{k}^{(1/2-\delta)(2+\iota)}} \left[\sum_{i=1}^{k-1} \frac{1}{i} \left\|\sum_{t=k+1}^{n} \varepsilon_{t-i} A_t\right\|_{2+\iota}\right]^{2+\iota}$$
$$\le O(1) \sum_{\tilde{k}=g_n}^{n-1} \frac{1}{\tilde{k}^{(1/2-\delta)(2+\iota)}} \left[\sum_{i=1}^{k-1} \frac{1}{i} \left(\sum_{t=k+1}^{n} \frac{1}{t}\right)^{1/2}\right]^{2+\iota}$$

$$\leq O(1) \sum_{\tilde{k}=g_n}^{n-1} \frac{1}{\tilde{k}^{(1/2-2\delta)(2+\iota)}} \left\{ \sum_{i=1}^{k-1} \frac{1}{i} \left[ \sum_{t=k+1}^{n} \frac{1}{(n-t+1)^{2\delta}t} \right]^{1/2} \right\}^{2+\iota}$$

$$\leq O(1) \sum_{\tilde{k}=g_n}^{n-1} \frac{1}{\tilde{k}^{(1/2-2\delta)(2+\iota)}} \left( \sum_{i=1}^{k-1} \frac{\log^{1/2} n}{in^{\delta}} \right)^{2+\iota} = o(1),$$

where the next-to-last step uses Lemma A.3. Similarly, we have

$$P\left( \max_{g_n \leq \tilde{k} < n} \frac{1}{\tilde{k}^{1/2-\delta}} \left| \sum_{i=k}^{n-1} \frac{1}{i} \sum_{t=i+1}^{n} \varepsilon_{t-i} A_t \right| > \epsilon \right)$$

$$\leq O(1) \sum_{\tilde{k}=g_n}^{n-1} \frac{1}{\tilde{k}^{(1/2-\delta)(2+\iota)}} \left[ \sum_{i=k}^{n-1} \frac{1}{i} \left( \sum_{t=i+1}^{n} \frac{1}{t} \right)^{1/2} \right]^{2+\iota} = o(1),$$

where we have used $\tilde{k} \geq n - i \geq n - t + 1$. By the preceding two inequalities, we have

(7.7) $$P\left( \max_{g_n \leq \tilde{k} < n} \frac{1}{\tilde{k}^{1/2-\delta}} \left| \sum_{t=k+1}^{n} v_t A_t \right| > \epsilon \right) = o(1).$$

Similarly, we can show that

(7.8) $$P\left( \max_{g_n \leq \tilde{k} < n} \frac{1}{\tilde{k}^{1/2-\delta}} \left| \sum_{t=k+1}^{n} \varepsilon_t A_{1t} \right| > \epsilon \right) = o(1).$$

Finally, by (7.3) and (7.6)–(7.8), we can show that the second part of Assumption 3.3(ii) holds. The first part of Assumption 3.3(ii) can be similarly proved and, hence, the details are omitted. □

## APPENDIX: LEMMAS A.1–A.3

We state three lemmas here whose proofs can be found in [28].

LEMMA A.1. *If Assumptions* 3.1(i), 3.2(i) *and* 3.3(i) *hold, then for any* $\eta > 0$

(a) $\lim_{n \to \infty} P\left( \max_{g_n \leq k \leq n} \sup_{|\lambda - \lambda_0| \geq \eta} \sum_{t=1}^{k} [\tilde{l}_t(\lambda) - \tilde{l}_t(\lambda_0)] + \sqrt{k} > 0 \right) = 0,$

(b) $\lim_{n \to \infty} P\left( \max_{g_n \leq n-k < n} \sup_{|\lambda - \lambda_0| \geq \eta} \sum_{t=k+1}^{n} [\tilde{l}_t(\lambda) - \tilde{l}_t(\lambda_0)] + \sqrt{n-k} > 0 \right) = 0.$

LEMMA A.2. *If the assumptions of Theorem* 3.1 *hold, then there exists a $\delta > 0$ such that*

(a) $\quad \max_{g_n \leq k \leq n} \dfrac{k^\delta}{n} \left| \sum_{t=1}^{k} \left[ \widetilde{D}_t(\hat{\lambda}_n(k)) \widetilde{D}'_t(\hat{\lambda}_n(k)) - \Omega \right] \right| = o_p(1),$

(b) $\quad \max_{g_n \leq n-k < n} \dfrac{(n-k)^\delta}{n} \left| \sum_{t=k+1}^{n} [\widetilde{D}_t(\hat{\lambda}_{1n}(k)) \widetilde{D}'_t(\hat{\lambda}_{1n}(k)) - \Omega] \right| = o_p(1).$

LEMMA A.3. *For any $u \in (0,1)$ and $v \in (0,\infty)$, it follows that*

$$\sum_{t=r+1}^{j} \frac{1}{(j-t+1)^u t^v} = O(1) \begin{cases} j^{1-u-v}, & \text{if } v < 1, \\ j^{-u} \log j, & \text{if } v = 1, \\ j^{-u} r^{1-v}, & \text{if } v > 1, \end{cases}$$

*where $O(1)$ holds uniformly in $j > r \geq 1$.*

**Acknowledgments.** The author greatly appreciates the very helpful comments of three referees, the Associate Editor and the Co-Editors M. L. Eaton and J. Fan.

DEPARTMENT OF MATHEMATICS
HONG KONG UNIVERSITY OF SCIENCE
AND TECHNOLOGY
CLEAR WATER BAY
HONG KONG
E-MAIL: maling@ust.hk